\documentclass[12pt]{amsart}
\usepackage{amssymb,amsmath,amsthm}

\theoremstyle{plain}
\newtheorem*{conjecture}{Conjecture}
\newtheorem*{thm}{Theorem}
\newtheorem{theorem}{Theorem}[section]
\newtheorem{proposition}[theorem]{Proposition}
\newtheorem{cor}[theorem]{Corollary}
\newtheorem{lemma}[theorem]{Lemma}

\theoremstyle{definition}
\newtheorem{definition}[theorem]{Definition}

\newtheorem{remark}[theorem]{\it Remark}

\newcommand{\dotminus}{\mathbin{\mathchoice%
{\buildrel .\lower.6ex\hbox{\vphantom{.}} \over {\smash-}}%
{\buildrel .\lower.6ex\hbox{\vphantom{.}} \over {\smash-}}%
{\buildrel .\lower.4ex\hbox{\vphantom{.}} \over {\smash-}}%
{\buildrel .\lower.3ex\hbox{\vphantom{.}} \over {\smash-}}}}%

\begin{document}

\title{Rational points on weighted projective spaces }  
\author{An-Wen Deng}

\address{ Institute of Mathematics\\
        Academia Sinica}
\email{adeng@math.sinica.edu.tw}

\date{}

\maketitle
\begin{abstract}
In this paper, we count the rational points on the weighted projective spaces defined over number fields w.r.t. `` size". An asymptotic formula which generalizes the result in \cite{Schn79} is obtained. Furthermore, we count also the rational points on the product of weighted projective spaces.
\end{abstract}
\section{Introduction}
This article is one of our papers discussing the rational points on the weighted projective varieties.  We use the argument in \cite{Schn79}  to count the number points on the weighted projective spaces over number fields.  As in our paper \cite{deng2}, we  recognized that the Weil  heights or Arakelov heights are not the  ``suitable"  functions to measure the ``size"  of the rational points on the weighted projective varieties because of weights.  Instead of heights, we introduced the function `` primitive size", $Size$,  on the weighted projective varieties . Here, we will  describe the function ``size"  more subtly. To each effective Weil divisor $D$ over the product of weighted projective spaces we can associate a counting  function ``$Size_D$".
\par
For a quasismooth weighted projective variety $X$  over the  number field $\Bbbk$ whose rational Picard group $Pic(X)\otimes{{\Bbb Q}}$ has rank one and its amplitude is $a$, we expect that 
$$
N(X(\Bbbk),T,Size^a)\sim CT
$$
where $N(X(\Bbbk),T,Size^a):=\#\{P\in X(\Bbbk):Size^a(P)<T\}$ and $C$ is a constant. This asymptotic formula is similar to the asymptotic formula in Manin's linear growth conjecture for the smooth Fano varieties having Picard group of rank one. 
The expected asymptotic formula coincides with our calculations if taking  $X={\bf P}(W)$  as a well-formed projective space (cf. Theorem (A)) or  as a quasismooth well-formed weighted complete intersection in \cite[Th. A]{deng2}. We remark here that the above mentioned quasismooth varieties are (singular) Fano and have at worst the quotient cyclic singularities. ( Cf. \cite{fletcher}, \cite{De75}, \cite{BR86} or  \cite{Do82})
\par  
Furthermore, we discuss also the rational points on the product of well-formed weighted projective spaces $ X=\prod_{i=1}^{\rho}{\bf P}(W_i)$. Since the product of weighted projective spaces has the rational Picard group of  rank $\rho$, the choice of ``sizes" is more complicated. W.r.t. ``sizes" associated to different divisorical  sheaves, we count the rational points and obtain the correspondent asymptotic formulas. Especially, when the ``anticanonical size" $Size_{-K}$ is chosen, the asymptotic formula (cf. Theorem (B)) becomes
$$
N(X(\Bbbk),T,Size_{-K})\sim CT(\log T)^{\rho-1}
$$
which is also related to the Manin conjecture \cite{FMT89}.
\section{Weighted projective spaces and their products}
In this section we introduce some notions about the weighted projective varieties and their products. For details about the weighted projective spaces cf. \cite{De75}, \cite{BR86}  \cite{Do82} and \cite{fletcher}.
\par
Let  $W=(w_1,\dots,w_m)$ be an $m$-tuple of positive integers. Let $F$ be an algebraically closed field of characteristic $0$. Denote $S$  as the polynomial ring $F[x_1,\dots,x_m]$. Regard $S=S(W)$ as a graded ring with the graduation ${\deg}(x_k):=w_k$ for $k=1,\dots,m$.   The {\it weighted projective space} ${\bf P}(W)$ w.r.t. weight $W$ is denoted by 
${\bf P}(W):=Proj(S(W))$.  The sum of  $w_i's$ is denoted by $|W|:=\sum_iw_i$.
If without explicite mentioning, the weighted projective space ${\bf P}(W)$ is always assumed to be {\it well-formed}, i.e. each $(m-1)$ elements from $W$ are relatively prime, because of  \cite[I.2.5. and I.2.7.]{fletcher}.
\par
Given weights  $W_1=\{w_{11},\cdots,w_{1,1_m}\}, \cdots,W_k=\{w_{k1},\cdots,w_{k,k_m}\}$, let 
${\bf P}(W_1)\times\cdots\times {\bf P}(W_k)
$ 
 be the product of weighted projective spaces.  It should be noted in more general context that the product of weighted projective spaces can be regarded as a toric variety. Here, we will not use this  treatment.  Let $S_i$ ($i=1,\dots,k$) be the polynomial ring $F[x_{i1},\cdots,x_{i,i_m}]$. The  associated graded polynomial ring is 
$$
S:=\bigotimes_{i=1}^kS_i (W_i)
$$ 
with the graduation 
$$
{\rm multideg}(x_{1,1_i}\otimes\cdots\otimes x_{k,k_i})=(w_{1,1_i},\cdots,w_{k,k_i}).
$$
More precisely,  decompose $S$ to 
$
S=\bigoplus_{d}S_d
$
where $d$ runs over the index set ${\Bbb Z}_{\geq 0}^k$ and $S_d$ denotes the set of weighted multihomogeneous polynomials of multidegree $d$.  Denote 
$$
S_+:=\bigoplus_{d\not=0}S_d.
$$
The product ${\bf P}(W_1)\times\cdots\times {\bf P}(W_k)$ is nothing else than the set of prime weighted multihomogeneous ideals in $S$ which do not contain $S_+$.
\par
Let $Z$ denote the affine subscheme 
$$
\bigcup_{i=1}^k Spec(S_1)\otimes \cdots Spec(S_{i-1})\otimes \{0\}\otimes Spec(S_{i+1}) \cdots \otimes Spec(S_k)
$$ 
in $Spec (S)$. Denote $U:=Spec(S)\setminus Z$.  Let $\pi:Spec(S)\setminus Z\to {\bf P}(W_1)\times\cdots\times {\bf P}(W_k)$ be the canonical projection.  Let $X$ be a closed subvariety in ${\bf P}(W_1)\times\cdots\times {\bf P}(W_k)
$  generated by the weighted multihomogeneous polynomials $f_1,\cdots ,f_r$ .  Denote $X:=V_+(f_1,\cdots,f_r)$. The inverse image of $X$ under $\pi$ is denoted by $C(X)^*$ and called to be the associated {\it quasi affine cone}. The Zariski closure of $C(X)^*$, denoted by $C(X)$, is called to be the associated {\it affine cone}. If $C(X)^*$ is non-singular, $X$ is said to be {\it quasismooth}.  
\begin{remark}
The product of weighted projective spaces  $\prod_{i=1}^k{\bf P}(W_i)$ can be interpreted as the geometrical quotient space $U/{\Bbb G}_m^k$  ( remember  ${\Bbb G}_m=Spec(F[T,T^{-1}])$) where the group ${\Bbb G}_m^k$ acts morphically on $Spec(S)$ via 
$$
S=\bigotimes_{i=1}^k S_i\to \bigotimes_{i=1}^k (S_i\otimes F[T_i,T_i^{-1}]) 
\qquad
(x_{ij}\to x_{ij}\otimes T_i^{w_{ij}}).
$$
\par
Actually, the group ${\Bbb G}_m^k(F)=({F}^{\times})^k$ acts on $U(F)$ via 
$$
(\lambda_1,\cdots,\lambda_k)_*(x_1,\cdots,x_k)=((\lambda_1)_*x_1,\cdots,(\lambda_k)_*x_k)
$$
where 
$$
(\lambda_i)_*x_i=(\lambda_i^{w_{i1}}x_{i1},\cdots,\lambda_i^{w_{i_m,1}}x_{i,1_m}).
$$
Each closed point in 
$\prod_{i=1}^k{\bf P}(W_i)(F)$ is exactly an orbit under the group action ${\Bbb G}_m^k(F)$ in $U(F)$.
\end{remark}
\begin{remark}
The product of weighted projective spaces is quasismooth have at most the cyclic quotient singularities.
\end{remark}
Denote by $K_X$ the canonical (Weil) divisor if  $X$ is a normal variety. Let ${\mathcal O}_X(D)$ denote the {\it divisorical sheaf} associated to the Weil divisor $D$ on $X$.
\begin{remark}
Let $X=\prod_{i=1}^k{\bf P}(W_i)$ be the product of weighted projective spaces .  Then the following  adjunction formula holds:
$$
{\mathcal O}_X(K_X)\cong {\mathcal O}_{{\bf P}(W_1)}(-|W_1|)\boxtimes \cdots \boxtimes {\mathcal O}_{{\bf P}(W_k)}(-|W_k|).
$$
Moreover, each Weil-divisor on $X$ is ${\Bbb Q}$-Cartier, and  $Cl(X)\cong {\Bbb Z}^k$ where each $k$-tuple $(a_1,\cdots,a_k)\in{{\Bbb  Z}}^k$ corresponds the divisorical sheaf 
$$
{\mathcal O}_{{\bf P}(W_1)}(a_1)\boxtimes\cdots\boxtimes {\mathcal O}_{{\bf P}(W_k)}(a_k),
$$
since $Cl({\bf P}(W_i))\cong {\Bbb Z}$ can be identified with the  group   generated by the isomorphic class $[{\mathcal O}_{\bf P}(W_i)(1)]$  \cite[prop. 2.3]{De75}.
\end{remark} 
\section{``Sizes'' on  the Product of Weighted Projective Spaces}
In rest of this paper we fix the following notations:
Let $\Bbbk$  denote   a number field,
  $O_{\Bbbk}$ the  ring of integers in  $\Bbbk$,  ${\frak U}(\Bbbk)$ the unit group, $ \mu(\Bbbk)$ the group of unity in $\Bbbk$, $D_{\Bbbk}$ the absolute value of the discriminant,  $V_{\Bbbk}$  the set of places over $\Bbbk$,  $S_{\infty}$ the set of archimedean places, $S_f$ the set of finite places,
$N:=[\Bbbk:{\Bbb Q}] $,
$N=r_1+2r_2 $   where $r_1$ denotes the number of real places and $r_2$ the number of complex places,  $R$ the regulator of $\Bbbk$, $Cl(\Bbbk)$ the group of ideal classes, $h$ the class number and $w:=\#\mu(\Bbbk)$.  Let the $v$-norms for $v\in S_{\Bbbk}$ be  normalized so that: (1) $|x|_v=|x|^{N_v}$  where $ N_v=1$  (resp. $2$) if $v$ is real (resp. complex) and (2) $|x|_v=|N_{k_v/{\Bbb Q}_p}(x)|_p$  as $v|p$ where $|p|_p=1/p$ . So, the product formula can be written as 
$$
 \prod_{v\in S_{\Bbbk}}|x|_v=1
$$
for $x\in \Bbbk^{\times}$. 
\par
Assume that the variety $X$ is in the product of  well formed  weighted projective space $\prod_{i=1}^k{\bf P}(W_i)$, where $W_i=(w_{i1},\dots, w_{i,i_m})$. 
Suppose  that ${\mathcal L}$  is an invertible sheaf over $X$. Then, according to the height machinery, one can define the Weil height or the  Arakelov height w.r.t. ${\mathcal L}$. Besides  heights, we introduce the other functions ``sizes'' for $\Bbbk$-rational points on the  varieties in the product of weighted projective spaces defined over $\Bbbk$.  We will see that the function ``size'' is the ``right'' function for counting  the rational points on  the  varieties in the product of weighted projective spaces  in contrast to  heights. 
 \begin{lemma}
The set of $\Bbbk$-rational points $\prod_{i=1}^k{\bf P}(W_i)(\Bbbk)$ can be regarded as the quotient space $\prod_{i=1}^k(\Bbbk^{i_m}\setminus\{0\})/({\Bbbk}^{\times})^k$, where $({\Bbbk}^{\times})^k$  acts on via $\prod_{i=1}^k(\Bbbk^{i_m}\setminus\{0\})$
$$
(a_1,\cdots,a_k)_*(x_1,\dots,x_k)=((a_1)_*x_1,\dots,(a_k)_*x_k). 
$$
In the other words, each $\Bbbk$-rational point $x=[x_1,\dots,x_k]$ on $\prod_{i=1}^k{\bf P}(W_i)$  is an orbit $Orb_x$ under the group action  $({\Bbbk}^{\times})^k$ .
\end{lemma}
\begin{proof}
 If $y=(y_1,\dots, y_k)\sim (x_1,\dots, x_k)=x  \in \prod_{i=1}^k(\Bbbk^{i_m}\setminus\{0\})$ then there exists an $a=(a_1,\cdots,a_k)\in (\overline{\Bbbk}^{\times})^k$ such that $a_*y=x$,  for each $j=1,\dots ,k$, $a_j$  may be assumed in  some finite field extension ${\Bbb L}/\Bbbk$.  Then, for some nonzero $a_i$, $\Bbbk a_i^{w_{ij}}$ are in $\Bbbk $. So, $\Bbbk a_i^{\gcd_j(w_{ij})}$ is in $\Bbbk $. But $\gcd_j(w_{ij})=1$ for ${\bf P}(W_i)$ being well-formed. This implies that $a\in \Bbbk $.
\end{proof}
The definition of  `` primitive size" for $x\in {\bf P}(W)(\Bbbk)$ is as follows.
\begin{definition}
For $x\in (\Bbbk)^m\setminus\{0\}$,
define  the  fractional ideal ${\frak I}(x)$  via ${\frak I}^{-1}(x):=\{a\in \Bbbk: a_*x\in {O}_{\Bbbk}\}$.
  Set $H_{\infty}(x):=\prod_{v\in S_{\infty}}\max_{i=1,\cdots ,m}|x_i|_v^{\frac{1}{w_i}}$.
Define 
$$
Size(x):=\frac{1}{N({\frak I}(x))}H_{\infty}(x).
$$
\end{definition}
\begin{proposition}
If  $x=(x_1,\cdots,x_m)\in O_{\Bbbk}^m$ then ${\frak I}(x)$ is really an integral idea of  $O_{\Bbbk}$.  
Moreover,  one has  the fact  ${\frak I}(a_*x)=a{\frak I}(x)$.  Consequently, both fractional ideals ${\frak I}(a_*x)$ and ${\frak I}(x)$ are in the same ideal class.
\end{proposition}
\begin{proof}
The first statement is clear from the definition of ${\frak I}(x)$. The second statement can be done from
\begin{eqnarray*}
b\in {\frak I}^{-1}(a_*x) \iff b_*(a_*x)\in O_{\Bbbk}^m
\iff (ba)_*x\in O_{\Bbbk}^m\\
 \iff ba \in {\frak I}^{-1}(x)\iff b\in a^{-1}{\frak I}^{-1}(x).
\end{eqnarray*}
\end{proof}
\begin{remark}
The well-definedness of the  function ``primitive size" is  a deduction from the product formula and the proposition above.
\end{remark}
\begin{remark}
For $x\in {\bf P}(W)(\Bbb Q)$, the ``primitive size'' of $x$ can be  given  by
$$
Size(x):=\min_{y\in Orb_x, y\in {\Bbb Z}^m}H_{\infty}(y).
$$
\par
The function ``primitive size'' coincides with the primitive Weil height when ${\bf P}(W)={\bf P}^{m-1}$ defined over $\Bbbk$, because in this case the ideal ${\frak I}(x)$ becomes the the ideal generated by the coordinates $x_1,\cdots, x_m$, where $\prod_{v\in S_f}\max_{i=1,\cdots,m}|x_i|=\frac{1}{N({\frak I}(x))}$.
\end{remark}
\begin{remark}
In general,  the function size is not equivalent to the  height.  For example, given a weighted projective space ${\bf P}(1,1,2)$ over ${\Bbb Q}$.  The twisted sheaf ${\mathcal O}(2)$ is very ample. Choose $P=(p,p,p)\in {\bf P}(1,1,2)({\Bbb Q})$ where $p$ is a prime number.  Then 
$$
H_{{\mathcal O}(2)}(P)=\frac{\max(|x_1|^2,|x_2|^2, |x_3|)}{\gcd(x_1^2,x_2^2,x_3)}\Bigg|_{(x_1,x_2,x_3)=(p,p,p)}=p
$$
and $Size^2(P)=p^2$.
\end{remark}
We define the size of rational points on  the product of weighted projective spaces $\prod_i{\bf P}(W_i)$ w.r.t. the Weil divisors. 
\begin{definition}
Let $D$ be an effective  Weil divisor  of  the product of weighted projective spaces $\prod_{i=1}^k{\bf P}(W_i)$ correspondent to $(a_1,\cdots,a_k)\in{\Bbb N}^r$. 
The function size w.r.t. the Weil divisor $D$ is defined via
$$
Size_{D}(x):=Size^{a_1}(x_1)\cdots Size^{a_k}(x_k)
$$
for $x=(x_1,\cdots,x_k)\in \prod_{i=1}^k{\bf P}(W_i)(\Bbbk)$.
\end{definition}
According to the adjunction formula, the anticanonical divisor $-K$ is correspondent to 
$( |W_1|,\cdots,|W_k|)$.  Therefore, $\prod_{i=1}^k{\bf P}(W_i)$ is Fano. ( Recall that a normal variety is {\it Fano } iff its anticanonical divisor is ${\Bbb Q}$-Cartier and ample.) The {\it anticanonical size} $Size_{-K}$ is given by
$$
Size_{-K}(x):=Size^{|W_1|}(x_1)\cdots Size^{ |W_k|}(x_k)
$$
for $x=(x_1,\cdots,x_k)\in \prod_{i=1}^k{\bf P}(W_i)(\Bbbk)$.
\section{Integral points modulo units}
We count the $\Bbbk$-rational points on the well-formed weighted projective space ${\bf P}(W)$ in following  three sections. 
\par
The  variety ${\bf A}^m\setminus \{0\}$ is the quasi  affine cone  of the weighted projective space ${\bf P}(W)$.  The subgroup ${\frak U}(\Bbbk)$ of $\Bbbk^{\times}$ acts on the set ${\bf A}^m(\Bbbk)\setminus \{0\}={\Bbbk}^m\setminus \{0\}$.  In this section, we generalize  the idea of \cite{Schn79} to  count  the integral points  on ${\bf A}^m\setminus\{0\}$ modulo ${\frak U}(\Bbbk)$ within a suitable weighted expanding compact set.
\par
For $c\in Cl({\Bbbk})$ an ideal class denote
$$
N(c,T):=\#\{x\in {\bf P}(W)(\Bbbk):Size(x)\leq T, [{\frak I}(x)]=c\}
$$
For brevity, choose  $\frak A$ as  an integral ideal in the ideal  $c$. One has
$$
{\frak A}^{-1}_*({\frak A}^{w_1}\times\cdots\times {\frak A}^{w_m})=O_{\Bbbk}^m.
$$
Consider  the number of orbits ${\frak U}(\Bbbk)_*(x)$   for $ x\in ({\frak A}^{w_1}\times\cdots\times {\frak A}^{w_m})\setminus \{0\}$ such that $H_{\infty}(x)\leq TN({\frak A})$. This number is denoted by  $N_{{\frak U}(\Bbbk)}({\frak A},T)$  and satisfies
$$
N_{{\frak U}(\Bbbk)}({\frak A},T)=\#\{ {\frak U}(\Bbbk)_*(x)\in \frac{{\Bbbk}^m-\{0\}}{{\frak U}(\Bbbk)}:{\frak I}(x)\subset {\frak A}, H_{\infty}(x)\leq TN({\frak A})\}.
$$
Denote 
$$
\overline{N}_{{\frak U}(\Bbbk)}({\frak A},T)=\#\{ {\frak U}(\Bbbk)_*(x)\in \frac{{\Bbbk}^m-\{0\}}{{\frak U}(\Bbbk)}:{\frak I}(x)={\frak A}, H_{\infty}(x)\leq TN({\frak A})\}.
$$
We note here that $N(c,T)=\overline{N}_{{\frak U}(\Bbbk)}({\frak A},T)$.  This can be done by the bijectivity of the mapping from $\{ {\frak U}(\Bbbk)_*(x)\in \frac{{\Bbbk}^m-\{0\}}{{\frak U}(\Bbbk)}:{\frak I}(x)={\frak A}, H_{\infty}(x)\leq TN({\frak A})\}$ to $\{x\in {\bf P}(W)(\Bbbk):Size(x)\leq T, [{\frak I}(x)]=c\}$ via ${\frak U}(\Bbbk)_*x\to \Bbbk_*^{\times}x$.  We count at first $N_{{\frak U}(\Bbbk)}({\frak A},T)$.
\par
Let $\Bbbk_{\infty}:=\Bbbk\otimes_{{\Bbb Q}}{\Bbb R}=\prod_{v\in S_{\infty}}\Bbbk_v.$
We recall here the Dirichlet's unit theorem.
\begin{thm}[Dirichlet]
Let the mapping
$\ell :{\frak U}(\Bbbk)\to {\Bbb R}^{r_1+r_2}$  be given via 
$u\to (\log |u|_v)_{v\in S_{\infty}}$. Then  $\Gamma:=Im(\ell)$ builds a lattice of rank $r:=r_1+r_2-1$ in the hyperplane $H$ defined by $\sum_{v\in S_{\infty}}y_v=0$ and $ker(\ell)=\mu(\Bbbk)$.
\end{thm}
Let $pr:{\Bbb R}^{r_1+r_2}\to H$ be the projection along the vector $(N_v)_{v\in S_{\infty}}$, where more precisely $pr_v(y)=y_v-(\frac{\sum_{w\in S_{\infty}}y_w}{N})N_v.$
Let 
$$
\eta:\prod_{v\in S_{\infty}}(\Bbbk_v^m\setminus \{0\})\to {\Bbb R}^{r_1+r_2}
$$ 
be given via
$\eta=(\eta_v)_{v\in S_{\infty}}$ where $\eta_v:\Bbbk_v^m\setminus \{0\}\to {\Bbb R}$  is defined by 
$$
\eta_v(Z_v)=\log\max_{i=1,\cdots,m}|Z_{v,i}|_v^{\frac{1}{w_i}}
$$
 for $Z=(Z_v)_{v\in S_{\infty}}=(Z_{v,i})_{v,i=1,\cdots,m} \in \prod_{v\in S_{\infty}}(\Bbbk^m_v\setminus \{0\}).$
Pick $F$ as a fundamental domain for $H/{\Gamma}$. Then $\Delta:=(pr\circ \eta)^{-1}F$ is a fundamental domain for $\prod_{v\in S_{\infty}}(\Bbbk^m \setminus\{0\})/{\frak U}(\Bbbk).$
We fix our choice of $F$, which is standard.
Let $\{u_1,\cdots, u_r\}$ be a basis for $\Gamma$ in $H$. Let $\{ \check{u}_1,\cdots, \check{u}_r\}$ be the dual basis, i.e. $\check{u}_j:H\to {\Bbb R}$   is  the linear functional satisfying 
$\check{u}_j(u_k)=\delta_{jk}.$ We choose $F=\{y\in H:0\leq \check{u}_j(y)< 1 \quad \forall j=1,\cdots,r\}.$
\par
Set the sets 
$$
{\frak D}(T):=\{Z\in \prod_{v\in S_{\infty}}(\Bbbk_v^m\setminus \{0\}):\prod_{v\in S_{\infty}}\max_{i=1,\cdots ,m}|Z_{v,i}|_v^{\frac{1}{w_i}}\leq T \}
$$
and ${\frak B}(T):={\frak D}(T)\cap \Delta$. We note here that ${\frak D}(T)$ is ${\frak U}(\Bbbk)$-stable. This can be done by
$$
\prod_v\max_i|u^{w_i}Z_{v,i}|_v^{\frac{1}{w_i}}=\prod_v|u|_v\prod_v\max_i|Z_{v,i}|_v^{\frac{1}{w_i}}=\prod_v\max_i|Z_{v,i}|^{\frac{1}{w_i}}.
$$
The properties concerning the fundamental domain and the group actions are as follows.
\begin{proposition}
(1) $ \Delta$ is $\mu(\Bbbk)$-stable.\\
(2) ${\frak U}(\Bbbk)_*(\Delta)=H.$\\
(3) $u_*(\Delta)\cap \Delta=\emptyset$ for $u\in {\frak U}(\Bbbk)\setminus \mu(\Bbbk)$.
\end{proposition}
\begin{proof}
This can be deduced from \cite[Lemma 1.]{Schn79}
\end{proof}
We want to count the number of ${\frak U}(\Bbbk)$-orbits in $({\frak A}^{w_1}\times\cdots\times {\frak A}^{w_m})\setminus \{0\}$ within the expanding set ${\frak D}(T)$. According to the proposition above and the fact of ${\frak D}(T)$ being ${\frak U}(\Bbbk)$-stable, it suffices to count the number of $\mu(\Bbbk)$-orbits in  $({\frak A}^{w_1}\times\cdots\times {\frak A}^{w_m})\setminus \{0\}\cap {\frak B}(T)$.  Since the group $\mu(\Bbbk)$ acts effectively,  each $\mu(\Bbbk)$ contains $w$ points.  Therefore
\begin{proposition} \label{e1}
$$
w\#\{ {\frak U}(\Bbbk)_*(Z)\in \frac{({\frak A}^{w_1}\times\cdots\times {\frak A}^{w_m})\setminus \{0\}}{{\frak U}(\Bbbk)}: \prod_{v\in S_{\infty}}\max_{i=1,\cdots,m}|Z_{v,i}|^{\frac{1}{w_i}}\leq T\}$$
is the number of points of ${\frak A}^{w_1}\times\cdots\times {\frak A}^{w_m}\setminus \{0\}$ in ${\frak B}(T)$.
\end{proposition}
\begin{lemma}
(1) $t_*\Delta=\Delta$ for $t\in {\Bbb R}^{\times}$\\
(2) ${\frak B}(T)=T^{\frac{1}{N}}_*{\frak B}(1)$
\end{lemma}
\begin{proof}
(Cf. \cite[Lemma 3. p. 438]{Schn79}) 
(1) For $Z\in \prod_v(\Bbbk_v^m\setminus \{0\})$ and  $t\in {\Bbb R}^{\times}$, we have 
$$\eta(t_*Z)=\log|t|(N_v)_{v\in S_{\infty}}+\eta(Z).$$
Since $\eta $ is linear and annihilates the vector $(N_v)_{v\in S_{\infty}}$, we have  $pr\circ \eta(t_*Z)=pr\circ \eta(Z)$. The fact (1) is done.\par
(2) Let  $\rho(Z)=\prod_{v\in S_{\infty}}\max_{i=1,\cdots,m}|Z_{v,i}|^{\frac{1}{w_i}}$. Then,  $$
\rho(t_*Z)=|t|^{\sum_v N_v}\rho(Z)=|t|^N\rho(Z).
$$ 
Using the fact (1), the fact (2) is achieved.
\end{proof}
Here, we need some notions from weighted expanding sets.   Let ${\frak P}$ be a bounded subset of ${\Bbb R}^k$ with $C^1$-{\it parametrizable boundary }(or {\em well-shaped boundary }), i.e. there exist finitely many $C^1$ mappings  from $[0,1]^{k-1}$ to $R^k$ such that the boundary $\partial {\frak P}$ is contained in the union of those images.  Given a weight $W=(w_1,\cdots,w_k)\in {\Bbb N}^k$. The weighted expanding set ${\frak P}(T)$ is defined by 
$$
{\frak P}(T):=T_*{\frak  P}:=\{(T^{w_1}x_1,\cdots,T^{w_k}x_k): (x_1,\cdots,x_k)\in {\frak P}\}.
$$  
We state here an observation.
\begin{proposition}\label{ee}
Notations as above.  Let $\Lambda $ be a lattice in ${\Bbb R}^k$. Then
\begin{eqnarray*}
\#\Lambda\cap {\frak P}(T)&=&\frac{Vol_k({\frak P}(T))}{\det \Lambda}+O(\frac{Vol_{k-1}(\partial {\frak P}(T))}{\det \Lambda})\\
&=&T^{|W|}\frac{Vol_k({\frak P})}{\det \Lambda}+O(T^{|W|-w_{\min}})
\end{eqnarray*}
where $w_{\min}:=\min_{i=1,\cdots,m}w_i$.
\end{proposition}
\begin{proof}
Let $D$ be  a  paralleltope fundamental domain of ${\Bbb R}^k/\Lambda$.  The paralleltopes $c+D$ ( $c\in \Lambda\cap {\frak P}(T))$ cover ${\frak P}(T)$. Moreover, there exists at most number  of $O(\frac{Vol_{k-1}(\partial {\frak P}(T))}{\det \Lambda})$  such paralleltopes which have nonempty intersection with $\partial {\frak P}$ and are not contained in ${\frak P}(T)$. 
Since $\frak P$ has the $C^1$-parametrizable boundary, it can be shown that  $\frac{Vol_{k-1}(\partial {\frak P}(T))}{\det \Lambda}=O(T^{|W|-w_{\min}})$.    Then it is not difficult to obtain the conclusion.
\end{proof}
 We apply this observation to count the number of ${\frak A}^{w_1}\times\cdots\times {\frak A}^{w_m}\setminus \{0\}$ in ${\frak B}(T)$ by
taking  ${\frak P}={\frak B}(1)\subset {\Bbb R}^{Nm}$, ${\frak P}(T)=T^{\frac{1}{N}}_* {\frak B}(1)$, $\Lambda=({\frak A}^{w_1}\times\cdots\times {\frak A}^{w_m})$. Consider  ${\frak A}\hookrightarrow \Bbbk_{\infty}=\prod_{v\in S_{\infty}}\Bbbk_v$ , where $\frak A$ build a lattice in $\Bbbk_{\infty}$ with determinant $\frac{\sqrt{D_{\Bbbk}}N({\frak A})}{2^{r_2}}$. The lattice $\Lambda $ has therefore  determinant 
$$\frac{D_{\Bbbk}^{\frac{m}{2}}N({\frak A})^{|W|}}{2^{mr_2}}.$$
In fact, we have by Proposition (\ref{e1})
\begin{proposition}
\begin{eqnarray*}
&\#&\{ {\frak U}(\Bbbk)_*(Z)\in \frac{({\frak A}^{w_1}\times\cdots\times {\frak A}^{w_m})\setminus \{0\}}{{\frak U}(\Bbbk)}: \prod_{v\in S_{\infty}}\max_{i=1,\cdots,m}|Z_{v,i}|^{\frac{1}{w_i}}\leq T\}
\\
&=& 
\frac{1}{w}(\frac{T}{N({\frak A})})^{|W|}\frac{2^{mr_2}}{D_{\Bbbk}^{\frac{m}{2}}}Vol_{Nm}({\frak B}(1))+O((\frac{T}{N({\frak A})})^{|W|-\frac{w_{\min}}{N}})
\end{eqnarray*}
\end{proposition}
\begin{cor}\label{as1}
$$
N_{{\frak U}(\Bbbk)}({\frak A},T)= 
\frac{1}{w}{T}^{|W|}\frac{2^{mr_2}}{D_{\Bbbk}^{\frac{m}{2}}}Vol_{Nm}({\frak B}(1))+O(T^{|W|-\frac{w_{\min}}{N}})
$$
\end{cor}
\begin{remark}
It remains to show that the bounded set ${\frak B}(1)$ has really the  $C^1$-parametrizable boundary. We will discuss about this in the next section.
\end{remark}
\section{Computations  on the volume of ${\frak B}(1)$}
In this section, we show that the bounded set ${\frak B}(1)$ has really the $C^1$-parametrizable boundary. Furthermore, we compute  the volume of ${\frak B}(1)$ . Recall
\begin{eqnarray*}
{\frak B}(1)&:=&{\frak B}_W(1)\\
&:=&\left\{Z\in \prod_{v\in S_{\infty}}(\Bbbk_v^m\setminus \{0\}):\begin{array}{cc}& \prod_{v\in S_{\infty}}
\max_{i=1,\cdots,m}|Z_{v,i}|_v^{\frac{1}{w_i}}\leq 1\\ 
&pr\circ \eta(Z)\in  F \end{array}
\right\}.
\end{eqnarray*}
Set
\begin{eqnarray*}
{\frak B}_0(1):=\left\{X\in \prod_{v\in S_{\infty}}(\Bbbk_v^m\setminus \{0\}):\begin{array}{cc}& \prod_{v\in S_{\infty}}
\max_{i=1,\cdots,m}|X_{v,i}|_v\leq 1\\ 
&pr\circ \overline{\eta}(Z)\in  F \end{array}
\right\}
\end{eqnarray*}
where $\overline{\eta}=(\overline{\eta}_v)_{v\in S_{\infty}}$ and $\overline{\eta}_v(X_v):=\log\max_{i=1,\cdots,m}|X_{v,i}|_v$.
\begin{lemma} \label{c1}
The bounded set ${\frak B}_0(1)$ has the $C^1$-parametrizable boundary.
\\
( Cf. \cite[Lemma 10]{Schn79})
\end{lemma}
\begin{proposition}
${\frak B}_W(1)$ has the $C^1$-parametrizable boundary.
\end{proposition}
\begin{proof}
Construct the mappings $\varphi_{\epsilon}: {\frak B}_0(1)\to{\frak B}_W(1)$ via $(X_{v,i})_{v,i}\to (\epsilon(v,i)X_{v,i}^{w_i})_{v,i}$ where $ \epsilon:S_{\infty}\times \{1,\cdots,m\}\to \{\pm 1\}.$  It is clear that the union $\bigcup_{\epsilon}Im \varphi_{\epsilon}$ covers ${\frak B}_W(1)$. From Lemma (\ref{c1}), ${\frak B}_0(1)$ has the $C^1$-parametrizable boundary, i.e. there exist finitely many $C^1$-mappings $\phi_i:[0,1]^{mN-1}\to {\Bbb R}^{mN}$ such that $\bigcup_i Im \phi_i\supset \partial {\frak B}_0(1).$ Therefore, $\bigcup_{\epsilon,i}\varphi_{\epsilon}\circ \phi_i([0,1]^{mN-1})\supset \partial {\frak B}_W(1).$
\end{proof}
We generalize the computations in  \cite[Proposition. p.443]{Schn79}  as follows:
\begin{proposition}\label{as2}
$$
Vol_{mN}({\frak B}_W(1))=2^{mr_1}\pi^{mr_2}R|W|^{r_1+r_2-1}
$$
\end{proposition}
\begin{proof}
Consider the polar coordinates $(\rho_{v,i},\theta_{v,i})_{v\in S_{\infty},i=1,\cdots m}$ where (1) $\rho_{v,i}=|Z_{v,i}|, \theta=\arg Z_{v,i}$ if $v$ is complex and (2)  $\rho_{v,i}=|Z_{v,i}|, \theta={\rm sign} Z_{v,i}$ if $v$ is real .  Then,
$$
Vol_{mN}({\frak B}_W(1))=2^{mr_1}\int_{{\Bbb D}}\prod_{v :\mbox{\small complex}}\rho_{v,i}\prod_{v\in S_{\infty}}d\rho_{v,i}\prod_{v: \mbox{\small complex}}d\theta_{v,i}
$$
where the integral domain $\Bbb D$ is given via
$ \rho_{v,i}\geq 0, 0\leq \theta_{v,i}< 2\pi$
satisfying 
$\prod_{v \in S_{\infty}}\max_{i=1,\cdots ,m}\rho_{v,i}^{\frac{N_v}{w_i}}\leq 1$ and $pr \circ \log \max_{i=1,\cdots m} \rho_{v,i}^{\frac{N_v}{w_i}}\in \overline{F}$.
\par
Integrating w.r.t. $\prod_{v: \mbox{\small complex}}d\theta_{v,i}$, we have  
$$
 2^{mr_1}(2\pi)^{mr_2}\int_{{\Bbb D}}\prod_{v,i}\rho_{v,i}^{N_v-1}d\rho_{v,i}.
$$
Decompose the integral domain as follows: ${\Bbb D}=\bigcup_{\sigma:S_{\infty}\to \{1,\cdots,m\}}{\Bbb D}_{\sigma}$
where ${\Bbb D}_{\sigma}$ is the subset of ${\Bbb D}$  such that $\max_{i=1,\cdots,m}\rho_{v,i}^{\frac{1}{w_i}}=\rho_{v,\sigma(v)}^{\frac{1}{w_{\sigma(v)}}}$ for $v\in S_{\infty}$. The different ${\Bbb D}_{\sigma}'s$ meet each other only on the lower dimensional real loci defined via the equations $\rho_{v,i}^{\frac{1}{w_i}}=\rho_{v,j}^{\frac{1}{w_j}}.$
\par
It suffices to compute $\int_{{\Bbb D}_{\sigma}}\prod_{v,i}\rho_{v,i}^{N_v-1}d\rho_{v,i}$.
Let $t_{v,i}=\rho_{v,i}^{N_v}$, so $dt_{v,i}=N_v\rho_{v,i}^{N_v-1}d\rho_{v,i}$. Therefore, $\int_{{\Bbb D}_{\sigma}}\prod_{v,i}\rho_{v,i}^{N_v-1}d\rho_{v,i}=\frac{1}{2^{mr_2}}\int \prod_{v,i}dt_{v,i}$. Integrating w.r.t. $dt_{v,i}$ for all $i\not=\sigma(v)$, 
$\int_0^{t_{v,\sigma(v)}^{\frac{w_i}{w_{\sigma(v)}}}}dt_{v,i}=t_{v,\sigma(v)}^{\frac{w_i}{w_{\sigma(v)}}}, $ then we have
$$
2^{mr_1}(2\pi)^{mr_2}\int_{{\Bbb D}_{\sigma}}\prod_{v,i}\rho_{v,i}^{N_v-1}d\rho_{v,i}=2^{mr_1}\pi^{mr_2}\int \prod_{v\in S_{\infty}}t_{v,\sigma(v)}^{\frac{|W|}{w_{\sigma(v)}}-1}dt_{v,\sigma(v)}.
$$
Substituting  $\tau_{v}=t_{v,\sigma(v)}^{\frac{1}{w_{\tau(v)}}}$ in the integral above, one has
$$
2^{mr_1}(2\pi)^{mr_2}\int_{{\Bbb D}_{\sigma}}\prod_{v,i}\rho_{v,i}^{N_v-1}d\rho_{v,i}=2^{mr_1}\pi^{mr_2}\int \prod_{v\in S_{\infty}}(w_{\sigma(v)}\tau_v^{|W|-1}d\tau_v)
$$
where the integral domain is  
$$
\{\tau=(\tau_v)\in {\Bbb R}_{\geq 0}^{r_1+r_2}|0 \leq \prod_{v\in S_{\infty}} \tau_v  \leq 1,0\leq \check{u}_j (pr \circ \log \tau_v)  \leq 1 \}.
$$
Via the change of variables, $u=\prod_{v\in S_{\infty}}\tau_v$ and $\xi_j=\check{u}_j(pr \circ \log \tau_v)$ for $j=1,\cdots,r$, where the Jacobian is $\pm R$ and the integral domain becomes $0\leq u,\xi_j\leq 1$. We have 
$$
\int_{{\Bbb D}_{\sigma}}2^{mr_1}(2\pi)^{mr_2}\prod_{v,i}\rho_{v,i}^{N_v-1}d\rho_{v,i}=2^{mr_1}\pi^{mr_2}R\frac{\prod_{v\in S_{\infty}}w_{\sigma(v)}}{|W|}.
$$
Summing up the integrals over all ${\Bbb D}_{\sigma}'s$ for  $\sigma:S_{\infty}\to \{1,\cdots,m\}$, this leads to the conclusion.
\end{proof}
\section{The Asymptotic Formula on ${\bf P}(W)(\Bbbk)$}
In this section we  obtain the asymptotic formula on the set of rational points ${\bf P}(W)(\Bbbk)$ w.r.t. `` primitive size". Let  $B$ be a countable set, and  $\phi$ a counting function, i.e. for every positive number $T$ the cardinality of the set of $x\in B$ with $\phi (x)\leq T$ is finite. Let $N(B,T,\phi)$ denote the cardinality of this finite set. We are interested especially in $N({\bf P}(W)(\Bbbk),T,Size)=\#\{x\in {\bf P}(W)(\Bbbk): Size(x)\leq T\}$ as $T\to \infty$.
\begin{thm}[A]
$$
N({\bf P}(W)(\Bbbk),T,Size)\sim   C_{\Bbbk}^WT^{|W|}
$$
where
 $$
C_{\Bbbk}^W:=\frac{h}{\zeta_{\Bbbk}(|W|)}(\frac{2^{r_1+r_2}\pi^{r_2}}{\sqrt{D_{\Bbbk}}})^m\frac{R}{w}|W|^{r_1+r_2-1}
$$
 and  the error term is 
$O(T^{|W|-\frac{w_{\min}}{N}})$ (if  ${\bf P}(W)(\Bbbk)={\bf P}^1({\Bbb Q})$  the error term is taken to be  $O(T\log T)$).
\end{thm}
\begin{proof}
From the discussions above, we have
$
\#\{x\in {\bf P}(W)(\Bbbk): Size(x)\leq T\}=\sum_{c\in Cl(\Bbbk)}N(c,T)=\sum_{[{\frak A}]\in Cl(\Bbbk),{\frak A} \mbox{ ideal}}\overline{N}_{{\frak U}(\Bbbk)}({\frak A},T).
$
Since 
$$
{N}_{{\frak U}(\Bbbk)}({\frak A},T)=\sum_{{\frak B}: \mbox{ideal}}\overline{N}_{{\frak U}(\Bbbk)}({\frak A}{\frak B},\frac{T}{N({\frak B})}),
$$
By the M\"obius inversion formula,
$$
\overline{N}_{{\frak U}(\Bbbk)}({\frak A},T)=\sum_{{\frak B}: \mbox{ideal}}\mu({\frak B}){N}_{{\frak U}(\Bbbk)}({\frak A}{\frak B},\frac{T}{N({\frak B})}),
$$
where $\mu(\cdot)$ is the M\"obius function on the monoid of the ideals in $O_{\Bbbk}$.
By Corollary (\ref{as1}) and Proposition (\ref{as2}), we have
$$
{N}_{{\frak U}(\Bbbk)}({\frak A}{\frak B},\frac{T}{N({\frak B})})\sim (\frac{2^{r_1+r_2}\pi^{r_2}}{\sqrt{D_{\Bbbk}}})^m\frac{R}{w}|W|^{r_1+r_2-1}(\frac{T}{N({\frak B})})^{|W|}.
$$
Using the fact $ \zeta_{\Bbbk}(s)^{-1}=\sum_{{\frak B}: \mbox{ideal}}\frac{\mu({\frak B})}{N({\frak B})^s}$ for $Re(s)>1$, we can obtain the main term for $\overline{N}_{{\frak U}(\Bbbk)}({\frak A},T)$, i.e. $ \frac{1}{\zeta_{\Bbbk}(|W|)}(\frac{2^{r_1+r_2}\pi^{r_2}}{\sqrt{D_{\Bbbk}}})^m\frac{R}{w}|W|^{r_1+r_2-1}T^{|W|}.$ The estimate of the error term is standard and cf. \cite[Lemma 12.]{Schn79}.  Summing up $\overline{N}_{{\frak U}(\Bbbk)}({\frak A},T)$ for $[{\frak A}]\in Cl(\Bbbk)$, we obtain the conclusion.
\end{proof}
\begin{cor}
$$
N({\bf P}(W)(\Bbbk),T,Size_{-K}(x)\leq T\}\sim C_{\Bbbk}^WT
$$
and 
$$
N({\bf P}(W)(\Bbbk),T,Size^e(x)\leq T\}\sim C_{\Bbbk}^WT^{\frac{|W|}{e}}
$$
for $e\in {\Bbb N}$.
\end{cor}
\begin{proof}
By the adjunction formula, ${\mathcal O}(-K)\cong {\mathcal O}(|W|)$,  $|W|$ is the amplitude of  the Fano variety ${\bf P}(W)$, therefore $Size_{-K}=Size^{|W|}$.
\end{proof}
\begin{remark}\label{acc}
In fact, for any  Zariski open $U$ in ${\bf P}(W)$ over $\Bbbk$, we have
$$
N(U(\Bbbk),T,Size)\sim   C_{\Bbbk}^WT^{|W|}.
$$
The proof of this fact is based on the same argument to prove Theorem (A) with only slight modifications.  Let $\pi\to {\bf A}^m\setminus \{0\} \to {\bf P}(W)$ be the canonical projection. Let $S:={\bf P}(W)\setminus U$.
It suffices to replace ${\frak D}(T)$ (resp. ${\frak B}_W(T))$ by ${\frak D}(T)\cap U(\Bbbk_{\infty})$ (resp. ${\frak B}_W(T)\cap U(\Bbbk_{\infty})$ ), where  ${\frak B}_W(1)\cap U(\Bbbk_{\infty})$ has the $C^1$-parametrizable boundary and the same volume as ${\frak B}_W(1)$.
\end{remark}
A variety $X$ over $\Bbbk$ has a Zariski closed subset $S$ as a {\em subvariety of accumulation points} w.r.t. the counting function $\phi$, iff
$$
\limsup \frac{\log N(X(\Bbbk),T,\phi)}{\log T}>\limsup \frac{\log N((X\setminus S)(\Bbbk),T,\phi)}{\log T}.
$$
\begin{proposition}
${\bf P}(W)$ has no subvariety of accumulation points w.r.t. $``Size"$.
\end{proposition}
\begin{proof}
Use Remark (\ref{acc}).
\end{proof}
\section{Rational Points on the Product of weighted Projective Spaces}
Based on the asymptotic formula in Theorem (A), we are able to  count the rational points on the product of weighted projective spaces w.r.t. different ``sizes". Our treatment is quite elementary which is based on the following fact.
\begin{lemma}
Let $X$  and $Y$ be countable sets , and $\phi_1$ (resp. $\phi_2$) the counting functions on $X$ ( resp. $Y$). Let $\phi=\phi_1\times \phi_2$ be a counting function on $X\times Y$. Assume that 
$$ 
N(X,T,\phi_1)=C_XT^{\alpha(X)}(\log T)^{\beta}+\Bigg\{\begin{array}{ll} O(T^{\alpha(X)}(\log T)^ {\beta-1})&  \mbox{ if }  \beta\in {\Bbb N}\\
O(T^{\alpha(X)-\delta_1}) &  \mbox{ if } \beta=0 \end{array} $$
and $N(Y,T,\phi_2)=C_YT^{\alpha(Y)}+O(T^{\alpha(Y)-\delta_2})$ where $\beta\geq 0$, $\delta_1>0$, $\delta_2>0$ and $\alpha(X)\geq \alpha(Y)$.
Then
 \begin{eqnarray*}
\begin{array}{l}
N(X\times Y,T, \phi)= \\
\left\{ \begin{array}{ll}\frac{C_XC_Y\alpha(X)}{\beta+1}T^{\alpha(X)}(\log T)^{\beta +1}+O(T^{\alpha(X)}(\log T)^{\beta}) & \mbox{ if } \alpha(X)=\alpha(Y)  \\
\frac{C_XC_Y\alpha(Y)}{\alpha(X)-\alpha(Y)}T^{\alpha(X)}(\log T)^{\beta}+O(T^{\alpha(X)}(\log T)^{\beta-1}) & \mbox{ if } \alpha(X)>\alpha(Y),\beta>0\\
\frac{C_XC_Y\alpha(Y)}{\alpha(X)-\alpha(Y)}T^{\alpha(X)}+O(T^{\alpha(X)-\delta}) & \mbox{ if } \alpha(X)>\alpha(Y),\beta=0\\
\end{array} \right.
\end{array}
\end{eqnarray*}
for some $\delta>0$.
\end{lemma}
\begin{proof}
The treatment is similar to \cite[ 2.Proposition.]{FMT89} and elementary.
\begin{eqnarray*}
N(X\times Y,T, \phi)&=&\sum_{x\in X}N(Y,\frac{T}{\phi_1(x)},\phi_2)\\
&=&\sum_{x\in X}C_Y(\frac{T}{\phi_1 (x)})^{\alpha(Y)}+O((\frac{T}{\phi_1(x)})^{\alpha(Y)-\delta})
\end{eqnarray*}
The error term has  the similar structure as the main term,  we discuss here only about the main term.
Let $a(i)$  denote the cardinality of  the set $\{x\in X:i-1<\phi_1(x)\leq i\}$.
Then
$$
N(X\times Y,T, \phi)\sim C_Y\sum_{j=1}^Ta(j)(\frac{T}{j})^{\alpha(Y)}+O((\frac{T}{j})^{\alpha(Y)-\delta}),
$$
where 
$$
\sum_{i=1}^Ma(i)=N(X,M,\phi_1)=C_XM^{\alpha(X)}(\log M)^{\beta}+O(M^{\alpha(X)}(\log M)^{\beta-1}).
$$
Using the Abel summation, we have
$$
N(X\times Y,T, \phi)\sim 
C_XC_YT^{\alpha(Y)}\sum_{j=1}^Tj^{\alpha(X)}(\log j)^{\beta}\int_j^{j+1}\alpha(Y)(\frac{1}{x})^{\alpha(Y)+1}dx.$$
Replacing the summation by the integral for the sake of  the mean value theorem, we obtain
\begin{eqnarray*}
N(X\times Y,T, \phi)\sim
C_XC_Y\alpha(Y)T^{\alpha(Y)}\int_1^{T+1}(\log x)^{\beta}x^{\alpha(X)-\alpha(Y)-1}dx \\
\sim \left\{ \begin{array}{ll}\frac{C_XC_Y\alpha(X)}{\beta+1}T^{\alpha(X)}(\log T)^{\beta +1} & \mbox{ if } \alpha(X)=\alpha(Y) \\
\frac{C_XC_Y\alpha(Y)}{\alpha(X)-\alpha(Y)}T^{\alpha(X)}(\log T)^{\beta} & \mbox{ if } \alpha(X)>\alpha(Y). \end{array} \right.
\end{eqnarray*}
The estimate for the error term is silimar to the estimate for the main term.
\end{proof}
Using this lemma iteratedly, we can give  the asymptotic formula  for $N(\prod_{i=1}^k{\bf P}(W_i)(\Bbbk),T,Size_D)$ where $D$ is  the Weil divisor corresponding to $(a_1,\cdots,a_k)\in {\Bbb N}^k$, especially:
\begin{thm}[B]
Let $P(W_i)=P(w_{i,1},\cdots ,w_{i,m})$  for $i=1,\dots ,k$. Then their product satisfies the asymptotic formula
$$
N(\prod_{i=1}^k{\bf P}(W_i)(\Bbbk),T,Size_{-K})\sim C(W_1,\cdots,W_k,\Bbbk)T(\log T)^{k-1}
$$
where 
\begin{eqnarray*}
& &C(W_1,\cdots,W_k,\Bbbk)=\\
&  &\quad\frac{h^k}{k!\prod_{i=1}^k\zeta_{\Bbbk}(|W_i|)}(\frac{2^{r_1+r_2}\pi^{r_2}}{\sqrt{D_{\Bbbk}}})^{\sum_{i=1}^ki_m}(\frac{R}{w})^k\prod_{i=1}^k|W_i|^{r_1+r_2-1}.
\end{eqnarray*}
\end{thm}
\begin{remark}
The product of weighted projective space $\prod_{i=1}^k {\bf P}(W_i)$  has no subvariety of accumulation points w.r.t. ``$Size_{-K}$''.
\end{remark}
For comparison with  our expectation, we formulate a version of Manin's linear growth conjecture:
\begin{conjecture}[Manin]
Let $X$ be a nonsingular projective  Fano variety defined over the number field $\Bbbk$. Let $H_{-K}$ be a height function w.r.t. the anticanonical divisor. $X(\Bbbk)$ is supposed to be dense in $X$ in the sense of Zariski topology. Let $\rho:=rank Pic(X)$. Let $U$ be a Zariski open set in $X$ having no accumulation points w.r.t. $H_{-K}$. Then
$$
N(U(\Bbbk),T,H_{-K})\sim C(X,\Bbbk)T(\log T)^{\rho-1}
$$
where $C(X,\Bbbk)$ is a positive constant.
{\rm (Cf.\cite{BM90})}\\
\end{conjecture}
\begin{remark}
The product of weighted projective spaces $X=\prod_{i=1}^k {\bf P}(W_i)$  is a (singular) projective  Fano variety and has no accumulation points w.r.t. the anticanonical size  ``$Size_{-K}$".  It is also easy to see that $X(\Bbbk)$ is Zariski dense in $X$ for any number field $\Bbbk$. Furthermore, $Pic(X)\cong {\Bbb Z}^k$.  The asymptotic formula in Theorem (B)  is parallel to the expectation of Manin's conjecture.
\end{remark}

\bibliographystyle{alpha}

\end{document}